\documentclass[a4paper, final, 12pt]{article}
\usepackage{amsmath, amssymb, latexsym, amscd, amsthm,amsfonts,amstext}
\usepackage[mathscr]{eucal}
\usepackage{graphicx}
\usepackage{subfig}
\usepackage{float}
\usepackage{color}
\usepackage{hyperref}
\usepackage[utf8]{inputenc}
\usepackage[english]{babel}

\numberwithin{equation}{section}
\input{tcilatex}
\begin{document}

\title{The Mean Field Games System: Carleman Estimates, Lipschitz Stability
and Uniqueness }
\author{Michael V. Klibanov \and Department of Mathematics and Statistics,
University of North \and Carolina at Charlotte, Charlotte, NC 28223, USA
\and email: mklibanv@uncc.edu}
\date{}
\maketitle

\begin{abstract}
An overdetermination is introduced in an initial condition for the second
order mean field games system (MFGS). This makes the resulting problem close
to the classical ill-posed Cauchy problems for PDEs. Indeed, in such a
problem and overdetermination in boundary conditions usually takes place. A
Lipschitz stability estimate is obtained. This estimate implies uniqueness.
A new Carleman estimate is derived. The second estimate the author calls
\textquotedblleft quasi-Carleman estimate", since it contains two test
functions rather than a single one in conventional Carleman estimates. These
two estimates play the key role. Carleman estimates were not applied to the
MFGS prior to the recent work of Klibanov and Aveboukh in \emph{arXiv}:
2302.10709, 2023. Applications are discussed.
\end{abstract}



\section{Introduction}

\label{sec:1}

Prior to the recent work \cite{KA}, Carleman estimates were not applied to
the mean field games system (MFGS). The author is also unaware about
previously obtained stability estimates for the MFGS. Both the idea and the
methodology of this paper came from the field of Ill-Posed and Inverse
Problems, see, e.g. \cite{BK,BukhKlib}, \cite{Klib84}-\cite{KL} for some
representative publications of the author in this field. Indeed, we consider
an overdetermination in the initial condition for the Bellman equation,
which is one of two nonlinear parabolic equations forming the second order
MFGS. This overdetermination makes the considered problem close to the
classical ill-posed Cauchy problems for PDEs, which routinely use
overdetermined boundary conditions, see, e.g. \cite{Isakov,KT,KL,LRS}. We
derive here a new Carleman estimate as well as an unusual quasi-Carleman
estimate (see this section below about some details).   

That extra initial condition can be considered as a result of a measurement.
Measurements always contain errors. We obtain the Lipschitz stability
estimate for the solution of the MFGS with respect to an error in that
measurement as well as with respect to possible errors in the conventional
terminal and initial conditions for this system. Our Lipschitz stability
estimate immediately implies uniqueness of the solution of the MFGS with
that extra initial condition.

The mean field games theory (MFG) was first introduced by Lasry and Lions 
\cite{Lasry_Lions_2006_I,Lasry_Lions_2006_II,LS} and Huang, Caines and Malham%
\'{e} \cite{Huang_Caines_Malhame_2007,Huang_Malhame_Caines_2006}. In these
seminal publications the second order system of two coupled nonlinear
parabolic equations was derived, i.e. the MFGS was derived. The first
equation, for the value function, is the Bellman equation. The second
equation is the Fokker-Plank equation for the function describing the
density of players. These two equations are coupled and form the MFGS.

The MFG theory studies the behavior of infinitely many agents/players, who
are trying to optimize their own payoffs depending on the dynamics of the
whole set of players. It is assumed that players act in a rational way and
that the dynamics of each player is independent on the dynamics of any other
player. The MFG theory has various applications in finance, economics,
pedestrians flocking, interactions of electric vehicles, etc., see, e.g. 
\cite{A,Cou,PP,Trusov}.

We consider the MFGS in a bounded domain $\Omega \subset \mathbb{R}^{n}$
with the zero Neumann boundary condition, meaning that there is no flux
through the boundary of both the value function and the agents. In the
conventional setting, the terminal condition for the value function and the
initial condition for the density function are given, see, e.g. \cite{A}.
However, uniqueness results are rare in this case because this setting
expresses the Nash equilibrium, which is usually not unique. Uniqueness of
the solution of the MFGS is proven only under some additional constraints,
like, e.g. a monotonicity assumption \cite{Bardi_Fischer_2017}. Therefore,
if avoiding such a constraint, then it is necessary to impose an additional
condition to guarantee uniqueness of the resulting problem.

As stated above, our additional initial condition is a result of a
measurement. Therefore, it contains an error. Thus, it is desirable to
obtain an accuracy estimate for that unique solution with respect to that
error. Using the terminology of the theory of Ill-Posed and Inverse
Problems, we call this \textquotedblleft stability estimate". It was
proposed in \cite{KA} to use the terminal condition for the density function
as that additional condition. \emph{The first} Lipschitz stability estimate
for the MFGS was obtained in \cite{KA}. In particular, that estimate implies
uniqueness of the solution of the MFGS with that overdetermination. To
obtain that result, two new Carleman estimates were derived in \cite{KA}.

In this paper we replace the knowledge of \cite{KA} of the terminal
condition for the density function with the knowledge of the initial
condition for the value function. This knowledge might be obtained via, e.g.
conducting a poll among agents at the initial moment of time. Indeed,
suppose that at the moment of time $\left\{ t=0\right\} $ a player at an
arbitrary fixed position $\overline{x}\in \Omega $ has an approximate idea
about the value function $u\left( x,0\right) $ of neighboring players, who
are located in a small neighborhood of the point $\overline{x}.$ Assuming
that the function $u\left( x,0\right) $ is continuous, we obtain an
approximation for the number $u\left( \overline{x},0\right) $ as $u\left( 
\overline{x},0\right) =\lim_{x\rightarrow \overline{x}}u\left( x,0\right) .$
Assuming that each player has that approximate idea, we obtain that an
approximation for the function $u\left( x,0\right) $ is known for all points 
$x$ located in the domain of interest $\Omega $. Hence, a poll, being
conducted among players at $\left\{ t=0\right\} ,$ would provide one with an
approximation for the initial condition $u\left( x,0\right) $ of the value
function for all positions $x\in \Omega $ of the players. Since conducting a
poll is a serious effort, then we assume that polling takes place only once,
i.e. at a single moment of time $\left\{ t=0\right\} .$

Suppose now that, using the above mentioned conventional terminal and
initial conditions as well as the function $u\left( x,0\right) ,$ the unique
solution of the MFGS is found even prior the game process has actually
started. Then one would know \underline{in advance} approximations for the
optimal behaviors of both the value and the distribution functions of the
process. From this standpoint, the importance of our Lipschitz stability
estimate is that it provides an \emph{a priori} estimate of the accuracy of
the solution of the MFGS depending on ones estimate of the error in the
input data $u\left( x,0\right) .$ For example, margins of errors of various
polls are often published in the news.

To obtain the desired Lipschitz stability result, we derive a new Carleman
estimate as well as a quasi-Carleman estimate. The first one is different
from the ones of \cite{KA} because Carleman Weight Functions (CWFs) are
different. These differences lead to the differences in the proof, as
compared with \cite{KA}. We use the term \textquotedblleft quasi" for the
second estimate since two test functions are involved in it rather than the
traditional single test function. The CWF is the function, which is involved
as the weight function in the Carleman estimate, see, e.g. \cite{KL}.

\emph{For the first time}, a Carleman estimate was derived in the seminal
paper of Carleman \cite{Carl}. Since then these estimates have been used for
proofs of uniqueness and stability results for ill-posed Cauchy problems for
PDEs, see, e.g. \cite{H,Isakov,KT,KL,LRS,Yam}. The publication \cite%
{BukhKlib} is the first one, in which Carleman estimates were introduced in
the field of Coefficient Inverse Problems (CIPs). The goal of \cite{BukhKlib}
was to apply Carleman estimates for proofs of global uniqueness and
stability results for multidimensional CIPs. This is because only local
uniqueness theorems for multidimensional CIPs were known prior to \cite%
{BukhKlib}.  The framework of \cite{BukhKlib} has been broadly used since
then by many authors with the same goal. Since the current paper is not a
survey of the method of \cite{BukhKlib}, then we refer here only to a few
publications, out of many, in which the methodology of \cite{BukhKlib} is
used for this purpose \cite{BK,ImYam1,Isakov,Klib84,Klib92,KT,Ksurvey,KL,Yam}%
. Furthermore, a modification of the idea of \cite{BukhKlib} is also
currently explored for constructions of a number of versions of the
so-called convexification globally convergent numerical method for CIPs,
see, e.g. \cite{KL}.

The following are three main mathematical difficulties in the MFGS:

\begin{enumerate}
\item Time is running in two opposite directions in two nonlinear parabolic
equations forming the MFGS. Therefore, the conventional theory of parabolic
equations is inapplicable here.

\item The nonlinearity of each of these equations.

\item The possible presence of a non-local term in the interaction part in
the Bellman equation.
\end{enumerate}

\textbf{Remark 1.1}. \emph{We note that the minimal smoothness conditions
are of a low priority in the theory of Ill-Posed and Inverse Problems, see,
e.g. \cite{Nov1,Nov2}, \cite[Lemmata 2.2.1 and 2.2.3]{Rom2}. Analogously,
they have a low priority in this paper as well as in \cite{KA}. The reason
of this is that those problems usually are highly nonlinear ones, including
this paper and \cite{KA}. Hence, they are very challenging ones even without
the minimal smoothness assumptions. Thus, it makes sense to obtain valuable
results first without minimal smoothness conditions and then work on
incorporating such conditions on a later stage of research.}

All functions considered below are real valued ones. In section 2 we state
the problem, which we address here. Section 3 is devoted to the proofs of
the above mentioned new Carleman estimate as well as the quasi Carleman
estimate. In section 4 we prove the Lipschitz stability estimate and
uniqueness for the problem formulated in section 2. Everywhere below $\beta
=const.>0.$ 

\section{Statement of the Problem}

\label{sec:2}

Let $x\in \mathbb{R}^{n}$\ denotes the position $x$\ of an agent and $t\geq 0
$ denotes time. Let $T>0$ be a number and $\Omega \subset \mathbb{R}^{n}$ be
a bounded domain with a piecewise smooth boundary $\partial \Omega $. Denote 
\begin{equation*}
Q_{T}=\Omega \times \left( 0,T\right) ,S_{T}=\partial \Omega \times \left(
0,T\right) .
\end{equation*}%
For $(x,t)\in Q_{T},$ we denote $u(x,t)$\ the value function and $p(x,t)$\
the distribution of agents at the point $x$\ and at the moment of time $t.$
Consider the system of infinitely many identical agents when the dynamics of
each agent is governed by the stochastic differential equation with
reflection 
\begin{equation}
dX(t)=\varkappa (X(t))\alpha (t)dt+\sqrt{2\beta }dW_{t}-n(X(t))dl(t),
\label{2.0}
\end{equation}%
where $X(t)$ is the state of an agent, $W_{t}$ is a Wiener process, $l$ is a
reflection process, which acts only at the boundary $\partial \Omega ,$ $%
\alpha (t)$ is a control, see below about the function $\varkappa (X(t))$.
The solution of (\ref{2.0}) is a pair $(X(\cdot ),l(\cdot ))$, we refer to
the book \cite{Pilipenko} for some details about the solution of equation (%
\ref{2.0}). Suppose that each agent is triyng to maximize his/her payoff.
This payoff is:%
\begin{equation*}
\mathbb{E}\Bigg[h(X(T))+\int_{0}^{T}\bigg(\frac{\alpha ^{2}(t)}{2}+G\left(
x,t,\dint\limits_{\Omega }K\left( x,y\right) p\left( y,t\right) dy,p\left(
x,t\right) \right) dt\Bigg].
\end{equation*}%
Here $\mathbb{E}$ is the mathematical expectation corresponding to the
Wiener process in (\ref{2.0}). Here the function $G$ is the interaction
term, which we describe below.

It is well known that by fixing an initial distribution of agents, one would
obtain the MFGS with the initial condition imposed on $p\left( x,t\right) $
and the terminal condition imposed on $u\left( x,t\right) $. This MFGS is:

\begin{equation*}
u_{t}(x,t)+\beta \Delta u(x,t)+\frac{\varkappa ^{2}(x)}{2}(\nabla
u(x,t))^{2}+
\end{equation*}%
\begin{equation}
+G\left( x,t,\dint\limits_{\Omega }K\left( x,y\right) p\left( y,t\right)
dy,p\left( x,t\right) \right) =0,\text{ }\left( x,t\right) \in Q_{T},
\label{2.1}
\end{equation}%
\begin{equation}
p_{t}(x,t)-\beta \Delta p(x,t)+\nabla \cdot (\varkappa ^{2}(x)p(x,t)\nabla
u(x,t))=0,\text{ }\left( x,t\right) \in Q_{T}.  \label{2.2}
\end{equation}

Thus, the time is running in two different directions in equations (\ref{2.1}%
) and (\ref{2.2}). The coefficient $\varkappa ^{2}(x)$\ plays the role of%
\textit{\ }the\textit{\ }elasticity of the medium. As stated above, the
function $F$ is the interaction term. The nonlocal interaction is given by
the integral operator with its kernel $K$. The local interaction is due to
the fourth argument of $F$. Hence, we consider a general case when both
local and non-local interactions are included in the interaction function $G$%
.

We assume the zero Neumann boundary conditions%
\begin{equation}
\partial _{n}u(x,t)\mid _{S_{T}}=\partial _{n}p(x,t)\mid _{S_{T}}=0,
\label{2.3}
\end{equation}%
where $n\left( x\right) $ is the unit outward looking normal vector at $%
\partial \Omega .$ In the conventional setting of the MFGS one assumes the
knowledge of the terminal condition for the function $u\left( x,t\right) $
and the initial condition for the function $p\left( x,t\right) ,$%
\begin{equation}
u\left( x,T\right) =u_{T}\left( x\right) ,\text{ }p\left( x,0\right)
=p_{0}\left( x\right) .  \label{2.4}
\end{equation}%
We, however, impose the additional initial condition on the function $%
u\left( x,t\right) ,$%
\begin{equation}
u\left( x,0\right) =u_{0}\left( x\right) .  \label{2.5}
\end{equation}

The function $u_{0}\left( x\right) $ can be interpreted as a result of a
poll conducted among players at the initial moment of time $\left\{
t=0\right\} ,$ see section 1 for a more detailed discussion. 

\textbf{Stability Estimate Problem}. \emph{Assume that there exist two
triples of functions (\ref{2.4}), (\ref{2.5}) }$\left( u_{T}^{\left(
i\right) }\left( x\right) ,p_{0}^{\left( i\right) }\left( x\right)
,u_{0}^{\left( i\right) }\left( x\right) \right) $\emph{, }$i=1,2$\emph{.
Also, suppose that there exist two pairs of solutions of problem (\ref{2.1}%
)-(\ref{2.5}) corresponding to these two triples: }%
\begin{equation*}
\text{ }\left( u^{\left( i\right) }\left( x,t\right) ,p^{\left( i\right)
}\left( x,t\right) \right) \in H_{0}^{2}\left( Q_{T}\right) \times
H_{0}^{2}\left( Q_{T}\right) ,i=1,2.
\end{equation*}%
\emph{\ Estimate appropriate norms of differences }$\left( u^{\left(
2\right) }-u^{\left( 1\right) }\right) \left( x,t\right) $\emph{\ and }$%
\left( p^{\left( 2\right) }-p^{\left( 1\right) }\right) \left( x,t\right) $%
\emph{\ via certain norms of the following differences:} 
\begin{equation*}
\left( u_{T}^{\left( 2\right) }-u_{T}^{\left( 1\right) }\right) \left(
x\right) ,\left( p_{0}^{\left( 2\right) }-p_{0}^{\left( 1\right) }\right)
\left( x\right) ,\left( u_{0}^{\left( 2\right) }-u_{0}^{\left( 1\right)
}\right) \left( x\right) .
\end{equation*}

Due to (\ref{2.3}) introduce the subspace $H_{0}^{2}\left( Q_{T}\right) $ of
the space $H^{2}\left( Q_{T}\right) $ as%
\begin{equation}
H_{0}^{2}\left( Q_{T}\right) =\left\{ v\in H^{2}\left( Q_{T}\right)
:\partial _{n}v(x,t)\mid _{S_{T}}=0\right\} .  \label{2.300}
\end{equation}%
By the trace theorem 
\begin{equation*}
v\left( x,0\right) ,v\left( x,T\right) \in H^{1}\left( \Omega \right) ,\text{
}\forall v\in H_{0}^{2}\left( Q_{T}\right) .
\end{equation*}

\section{Carleman Estimates}

\label{sec:3}

As stated in section 1, both the Carleman estimate and the quasi Carleman
estimate of this section are new and their proofs are different from the
ones of \cite{KA}. It is well known that a Carleman estimate for a PDE
operator depends only on the principal part of that operator, see, e.g. \cite%
[Lemma 2.1.1]{KL}. Thus, we derive below a Carleman estimate for the
operator $\partial _{t}+\beta \Delta $ in equation (\ref{2.1}). As to
equation (\ref{2.2}), we derive a Carleman estimate for the operator $%
\partial _{t}-\beta \Delta $ being perturbed by a term depending on another
function.

Let $a>2$ be a number and $\lambda >2$ be a sufficiently large parameter. We
will choose parameters $a$ and $\lambda $ later. Introduce the Carleman
Weight Function $\varphi _{\lambda }\left( t\right) ,$%
\begin{equation}
\varphi _{\lambda }\left( t\right) =e^{2\left( T-t+a\right) ^{\lambda }},%
\text{ }t\in \left( 0,T\right) .  \label{3.1}
\end{equation}%
Hence, 
\begin{equation}
\min_{\left[ 0,T\right] }\varphi _{\lambda }\left( t\right) =\varphi
_{\lambda }\left( T\right) =e^{2a^{\lambda }},\text{ }\max_{\left[ 0,T\right]
}\varphi _{\lambda }\left( t\right) =\varphi _{\lambda }\left( 0\right)
=e^{2\left( T+a\right) ^{\lambda }}.  \label{3.2}
\end{equation}

\subsection{The Carleman estimate for the operator $\partial _{t}+\protect%
\beta \Delta $}

\label{sec:3.1}

\textbf{Theorem 3.1.}\emph{\ Let in (\ref{3.1}) the parameter }$a>2.$\emph{\
Choose the number }$\lambda _{0}$\emph{\ as }%
\begin{equation}
\lambda _{0}=\lambda _{0}\left( a\right) =16\left( T+a\right)
^{2}>16a^{2}>64.  \label{3.9}
\end{equation}%
\emph{Then\ the following Carleman estimate is valid  for all functions }$%
u\in H_{0}^{2}\left( Q_{T}\right) $\emph{:}%
\begin{equation*}
\dint\limits_{Q_{T}}\left( u_{t}+\beta \Delta u\right) ^{2}\varphi _{\lambda
}dxdt\geq 
\end{equation*}%
\begin{equation}
\geq \frac{2}{3}\sqrt{\lambda }\beta \dint\limits_{Q_{T}}\left( \nabla
u\right) ^{2}\varphi _{\lambda }dxdt+\frac{\lambda ^{2}}{12}a^{\lambda
-2}\dint\limits_{Q_{T}}u^{2}\varphi _{\lambda }dxdt+  \label{3.300}
\end{equation}%
\begin{equation*}
-\frac{2}{3}e^{2a^{\lambda }}\dint\limits_{\Omega }\left( \beta \left(
\nabla u\right) ^{2}+\frac{1}{2}u^{2}\right) \left( x,T\right) -\frac{2}{3}%
\lambda \left( T+a\right) ^{\lambda -1}\dint\limits_{\Omega }u^{2}\left(
x,0\right) dx,\text{ }\forall \lambda \geq \lambda _{0}.
\end{equation*}

\textbf{Proof. }Consider the function $v\left( x,t\right) ,$ 
\begin{equation}
v\left( x,t\right) =u\left( x,t\right) e^{\left( T-t+a\right) ^{\lambda }}
\label{3.3}
\end{equation}%
and express appropriate derivatives of the function $u$ through the
derivatives of the function $v$. We have: 
\begin{equation*}
u=ve^{-\left( T-t+a\right) ^{\lambda }},
\end{equation*}%
\begin{equation}
u_{t}=\left( v_{t}+\lambda \left( T-t+a\right) ^{\lambda -1}v\right)
e^{\left( T-t+a\right) ^{\lambda }},\text{ }\Delta u=\left( \Delta v\right)
e^{\left( T-t+a\right) ^{\lambda }}.  \label{3.4}
\end{equation}%
Hence,%
\begin{equation*}
\left( u_{t}+\beta \Delta u\right) ^{2}\varphi _{\lambda }=\left[
v_{t}+\left( \lambda \left( T-t+a\right) ^{\lambda -1}v+\beta \Delta
v\right) \right] ^{2}\geq 
\end{equation*}%
\begin{equation*}
\geq 2v_{t}\left( \lambda \left( T-t+a\right) ^{\lambda -1}v+\beta \Delta
v\right) =
\end{equation*}%
\begin{equation*}
=\left( \lambda \left( T-t+a\right) ^{\lambda -1}v^{2}\right) _{t}+\lambda
\left( \lambda -1\right) \left( T-t+a\right) ^{\lambda -2}v^{2}+
\end{equation*}%
\begin{equation*}
+\dsum\limits_{i=1}^{n}\left( 2\beta v_{t}v_{x_{i}}\right)
_{x_{i}}-\dsum\limits_{i=1}^{n}\dsum\limits_{i=1}^{n}\left( 2\beta
v_{x_{i}t}v_{x_{i}}\right) \geq 
\end{equation*}%
\begin{equation*}
\geq \frac{\lambda ^{2}}{2}\left( T-t+a\right) ^{\lambda -2}v^{2}+\left(
\lambda \left( T-t+a\right) ^{\lambda -1}v^{2}-\beta \left( \nabla v\right)
^{2}\right) _{t}+
\end{equation*}%
\begin{equation*}
+\dsum\limits_{i=1}^{n}\left( 2\beta v_{t}v_{x_{i}}\right) _{x_{i}}.
\end{equation*}%
We have used here the fact that $\lambda >2,$ which implies that $\lambda
-1>\lambda /2.$ Thus, 
\begin{equation*}
\left( u_{t}+\beta \Delta u\right) ^{2}\varphi _{\lambda }\geq \frac{\lambda
^{2}}{2}\left( T-t+a\right) ^{\lambda -2}v^{2}+
\end{equation*}%
\begin{equation}
+\left( \lambda \left( T-t+a\right) ^{\lambda -1}v^{2}-\beta \left( \nabla
v\right) ^{2}\right) _{t}+\dsum\limits_{i=1}^{n}\left( 2\beta
v_{t}v_{x_{i}}\right) _{x_{i}}.  \label{3.5}
\end{equation}%
Integrate (\ref{3.5}) over $Q_{T}$ using (\ref{3.1}) and (\ref{3.3}). Since $%
v\in H_{0}^{2}\left( Q_{T}\right) ,$ then Gauss formula and (\ref{2.300})
imply that the integral over $Q_{T}$ of the last term of (\ref{3.5}) equals
\ zero. We obtain 
\begin{equation*}
\dint\limits_{Q_{T}}\left( u_{t}+\beta \Delta u\right) ^{2}\varphi _{\lambda
}dxdt\geq \frac{\lambda ^{2}}{2}\dint\limits_{Q_{T}}\left( T-t+a\right)
^{\lambda -2}u^{2}\varphi _{\lambda }dxdt+
\end{equation*}%
\begin{equation}
-\beta e^{2a^{\lambda }}\dint\limits_{\Omega }\left( \nabla u\right)
^{2}\left( x,T\right) dx-\lambda \left( T+a\right) ^{\lambda
-1}\dint\limits_{\Omega }u^{2}\left( x,0\right) dx,\text{ }\forall \lambda
>2.  \label{3.6}
\end{equation}

We now need to incorporate in (\ref{3.6}) a non-negative term containing $%
\left( \nabla u\right) ^{2}.$ To do this, consider the following expression:%
\begin{equation*}
\left( -u_{t}-\beta \Delta u\right) ue^{2\left( T-t+a\right) ^{\lambda
}}=\left( -\frac{u^{2}}{2}e^{2\left( T-t+a\right) ^{\lambda }}\right)
_{t}-\lambda \left( T-t+a\right) ^{\lambda -1}u^{2}e^{2\left( T-t+a\right)
^{\lambda }}+
\end{equation*}%
\begin{equation*}
+\dsum\limits_{i=1}^{n}\left( -\beta uu_{x_{i}}e^{2\left( T-t+a\right)
^{\lambda }}\right) _{x_{i}}+\beta \left( \nabla u\right) ^{2}e^{2\left(
T-t+a\right) ^{\lambda }}.
\end{equation*}%
Integrating this identity over $Q_{T}$ and taking into account using (\ref%
{3.1}), we obtain%
\begin{equation*}
\dint\limits_{Q_{T}}\left( -u_{t}-\beta \Delta u\right) u\varphi _{\lambda
}dxdt=\beta \dint\limits_{Q_{T}}\left( \nabla u\right) ^{2}\varphi _{\lambda
}dxdt-\lambda \dint\limits_{Q_{T}}\left( T-t+a\right) ^{\lambda
-1}u^{2}\varphi _{\lambda }dxdt-
\end{equation*}%
\begin{equation}
-\frac{1}{2}e^{2a^{\lambda }}\dint\limits_{\Omega }u^{2}\left( x,T\right) dx+%
\frac{1}{2}e^{2\left( T+a\right) ^{\lambda }}\dint\limits_{\Omega
}u^{2}\left( x,0\right) dx.  \label{3.60}
\end{equation}%
Multiply identity (\ref{3.60}) by $\sqrt{\lambda }$ and sum up with (\ref%
{3.6}). We obtain%
\begin{equation*}
\sqrt{\lambda }\dint\limits_{Q_{T}}\left( -u_{t}-\beta \Delta u\right)
u\varphi _{\lambda }dxdt+\dint\limits_{Q_{T}}\left( u_{t}+\beta \Delta
u\right) ^{2}\varphi _{\lambda }dxdt\geq 
\end{equation*}%
\begin{equation*}
\geq \sqrt{\lambda }\beta \dint\limits_{Q_{T}}\left( \nabla u\right)
^{2}\varphi _{\lambda }dxdt+
\end{equation*}%
\begin{equation*}
+\frac{\lambda ^{2}}{2}\dint\limits_{Q_{T}}\left( T-t+a\right) ^{\lambda
-2}\left( 1-\frac{2\left( T-t+a\right) }{\sqrt{\lambda }}\right)
u^{2}\varphi _{\lambda }dxdt-
\end{equation*}%
\begin{equation}
-e^{2a^{\lambda }}\dint\limits_{\Omega }\left( \beta \left( \nabla u\right)
^{2}+\frac{1}{2}u^{2}\right) \left( x,T\right) -\lambda \left( T+a\right)
^{\lambda -1}\dint\limits_{\Omega }u^{2}\left( x,0\right) dx.  \label{3.7}
\end{equation}%
It follows from (\ref{3.9}) that 
\begin{equation*}
\frac{2\left( T-t+a\right) }{\sqrt{\lambda }}\leq \frac{2\left( T+a\right) }{%
\sqrt{\lambda }}\leq \frac{1}{2},\text{ }\forall \lambda \geq \lambda _{0}.
\end{equation*}%
Keeping in mind that $\left( T-t+a\right) ^{\lambda -2}\geq a^{\lambda -2},$
we obtain that (\ref{3.7}) leads to%
\begin{equation*}
\sqrt{\lambda }\dint\limits_{Q_{T}}\left( -u_{t}-\beta \Delta u\right)
u\varphi _{\lambda }dxdt+\dint\limits_{Q_{T}}\left( u_{t}+\beta \Delta
u\right) ^{2}\varphi _{\lambda }dxdt\geq 
\end{equation*}%
\begin{equation}
\geq \sqrt{\lambda }\beta \dint\limits_{Q_{T}}\left( \nabla u\right)
^{2}\varphi _{\lambda }dxdt+\frac{\lambda ^{2}}{4}a^{\lambda
-2}\dint\limits_{Q_{T}}u^{2}\varphi _{\lambda }dxdt-  \label{3.10}
\end{equation}%
\begin{equation*}
-e^{2a^{\lambda }}\dint\limits_{\Omega }\left( \beta \left( \nabla u\right)
^{2}+\frac{1}{2}u^{2}\right) \left( x,T\right) -\lambda \left( T+a\right)
^{\lambda -1}\dint\limits_{\Omega }u^{2}\left( x,0\right) dx,\text{ }\forall
\lambda \geq \lambda _{0}.
\end{equation*}%
Since $a>2$, then by (\ref{3.9}) 
\begin{equation}
\frac{\lambda }{2}\leq \frac{\lambda ^{2}}{8}a^{\lambda -2},\text{ }\forall
\lambda \geq \lambda _{0}>64.  \label{3.11}
\end{equation}%
By Cauchy-Schwarz inequality%
\begin{equation*}
\sqrt{\lambda }\dint\limits_{Q_{T}}\left( -u_{t}-\beta \Delta u\right)
u\varphi _{\lambda }dxdt+\dint\limits_{Q_{T}}\left( u_{t}+\beta \Delta
u\right) ^{2}\varphi _{\lambda }dxdt\leq 
\end{equation*}%
\begin{equation}
\leq \frac{\lambda }{2}\dint\limits_{Q_{T}}u^{2}\varphi _{\lambda }dxdt+%
\frac{3}{2}\dint\limits_{Q_{T}}\left( u_{t}+\beta \Delta u\right)
^{2}\varphi _{\lambda }dxdt.  \label{3.12}
\end{equation}%
Combining (\ref{3.12}) with (\ref{3.10}) and (\ref{3.11}), we obtain the
target estimate (\ref{3.300}) of this theorem. $\square $

\subsection{The quasi-Carleman estimate}

\label{sec:3.2}

\textbf{Theorem 3.2.} \emph{Let in (\ref{3.1}) the parameter }$a>2$\emph{\
and let the number }$\lambda _{0}$\emph{\ be the same as in (\ref{3.9}). Let
the function }$f\left( x,t\right) $\emph{\ be differentiable in }$\overline{Q%
}_{T}$\emph{\ with respect to }$x$\emph{\ and such that}%
\begin{equation}
\sup_{Q_{T}}\left\vert f\left( x,t\right) \right\vert ,\text{ }%
\sup_{Q_{T}}\left\vert \nabla f\left( x,t\right) \right\vert \leq C_{f}
\label{3.120}
\end{equation}%
\emph{for a certain number }$C_{f}>0.$\emph{\ Then the following
quasi-Carleman estimate holds for any two functions }$u,q\in H_{0}^{2}\left(
Q_{T}\right) $\emph{:}%
\begin{equation*}
\dint\limits_{Q_{T}}\left( u_{t}-\beta \Delta u+f\Delta q\right) ^{2}\varphi
_{\lambda }\geq \frac{\lambda ^{2}}{4}a^{2\lambda
-2}\dint\limits_{Q_{T}}u^{2}\varphi _{\lambda }dxdt+\beta \lambda a^{\lambda
-1}\dint\limits_{Q_{T}}\left( \nabla u\right) ^{2}\varphi _{\lambda }dxdt-
\end{equation*}%
\begin{equation}
-C_{1}\lambda \left( T+a\right) ^{\lambda -1}\dint\limits_{Q_{T}}\left(
\nabla q\right) ^{2}\varphi _{\lambda }dxdt-  \label{3.121}
\end{equation}%
\begin{equation*}
-\lambda \left( T+a\right) ^{\lambda -1}e^{2\left( T+a\right) ^{\lambda
}}\dint\limits_{\Omega }u^{2}\left( x,0\right) dx,\text{ }\forall \lambda
\geq \lambda _{0},
\end{equation*}%
\emph{where the number }$C_{1}=C_{1}\left( C_{f},\beta \right) >0$\emph{\
depends only on listed parameters.}

\textbf{Proof}. In this proof $C_{1}=C_{1}\left( C_{f},\beta \right) >0$
denotes different numbers depending only on listed parameters. Use again
substitution (\ref{3.3}). In addition, denote 
\begin{equation}
w\left( x,t\right) =q\left( x,t\right) e^{\left( T-t+a\right) ^{\lambda }}.
\label{3.13}
\end{equation}%
Obviously, 
\begin{equation}
v,w\in H_{0}^{2}\left( Q_{T}\right) .  \label{3.14}
\end{equation}%
We have%
\begin{equation*}
\left( u_{t}-\beta \Delta u+f\Delta q\right) ^{2}\varphi _{\lambda }=\left[
\lambda \left( T-t+a\right) ^{\lambda -1}v+\left( v_{t}-\beta \Delta
v+f\Delta w\right) \right] ^{2}\geq
\end{equation*}%
\begin{equation*}
\geq \lambda ^{2}\left( T-t+a\right) ^{2\lambda -2}v^{2}+2\lambda \left(
T-t+a\right) ^{\lambda -1}v\left( v_{t}-\beta \Delta v+f\Delta w\right) =
\end{equation*}%
\begin{equation*}
=\lambda ^{2}\left( T-t+a\right) ^{2\lambda -2}v^{2}+\left( \lambda \left(
T-t+a\right) ^{\lambda -1}v^{2}\right) _{t}-\lambda \left( \lambda -1\right)
\left( T-t+a\right) ^{\lambda -1}v^{2}+
\end{equation*}%
\begin{equation}
+\dsum\limits_{i=1}^{n}\left( -2\beta \lambda \left( T-t+a\right) ^{\lambda
-1}vv_{x_{i}}\right) _{x_{i}}+2\beta \lambda \left( T-t+a\right) ^{\lambda
-1}\left( \nabla v\right) ^{2}+  \label{3.15}
\end{equation}%
\begin{equation*}
+\dsum\limits_{i=1}^{n}\left( 2\lambda \left( T-t+a\right) ^{\lambda
-1}vfw_{x_{i}}\right) _{x_{i}}-
\end{equation*}%
\begin{equation*}
-2\lambda \left( T-t+a\right) ^{\lambda -1}f\nabla v\nabla w-2\lambda \left(
T-t+a\right) ^{\lambda -1}v\nabla f\nabla w.
\end{equation*}%
Since $a>2$ and $\lambda \geq \lambda _{0},$ then 
\begin{equation}
\left( \lambda ^{2}\left( T-t+a\right) ^{2\lambda -2}-\lambda \left( \lambda
-1\right) \left( T-t+a\right) ^{\lambda -1}\right) v^{2}\geq \frac{1}{2}%
\lambda ^{2}\left( T-t+a\right) ^{2\lambda -2}v^{2}.  \label{3.16}
\end{equation}%
Hence, using (\ref{3.15}) and (\ref{3.16}), we obtain%
\begin{equation*}
\left( u_{t}-\beta \Delta u+f\Delta q\right) ^{2}\varphi _{\lambda }\geq 
\frac{1}{2}\lambda ^{2}\left( T-t+a\right) ^{2\lambda -2}v^{2}+2\beta
\lambda \left( T-t+a\right) ^{\lambda -1}\left( \nabla v\right) ^{2}-
\end{equation*}%
\begin{equation*}
-2\lambda \left( T-t+a\right) ^{\lambda -1}f\nabla v\nabla w-2\lambda \left(
T-t+a\right) ^{\lambda -1}v\nabla f\nabla w
\end{equation*}%
\begin{equation}
+\left( \lambda \left( T-t+a\right) ^{\lambda -1}v^{2}\right) _{t}+
\label{3.17}
\end{equation}%
\begin{equation*}
+\dsum\limits_{i=1}^{n}\left( -2\beta \lambda \left( T-t+a\right) ^{\lambda
-1}vv_{x_{i}}+2\lambda \left( T-t+a\right) ^{\lambda -1}vfw_{x_{i}}\right)
_{x_{i}}.
\end{equation*}

Next, by Cauchy-Schwarz inequality%
\begin{equation*}
-2\lambda \left( T-t+a\right) ^{\lambda -1}f\nabla v\nabla w-2\lambda \left(
T-t+a\right) ^{\lambda -1}v\nabla f\nabla w\geq 
\end{equation*}%
\begin{equation}
\geq -\beta \lambda \left( T-t+a\right) ^{\lambda -1}\left( \nabla v\right)
^{2}-\lambda \frac{C_{f}^{2}}{\beta }\left( T-t+a\right) ^{\lambda -1}\left(
\nabla w\right) ^{2}-  \label{3.18}
\end{equation}%
\begin{equation*}
-\lambda \left( T-t+a\right) ^{\lambda -1}v^{2}-\lambda C_{f}^{2}\left(
T-t+a\right) ^{\lambda -1}\left( \nabla w\right) ^{2}.
\end{equation*}%
Hence, (\ref{3.3}), (\ref{3.13}), (\ref{3.17}) and (\ref{3.18}) imply%
\begin{equation*}
\left( u_{t}-\beta \Delta u+f\Delta p\right) ^{2}\varphi _{\lambda }\geq 
\frac{\lambda ^{2}}{4}a^{2\lambda -2}u^{2}\varphi _{\lambda }+\beta \lambda
a^{\lambda -1}\left( \nabla u\right) ^{2}\varphi _{\lambda }-
\end{equation*}%
\begin{equation}
-C_{1}\lambda \left( T+a\right) ^{\lambda -1}\left( \nabla q\right)
^{2}\varphi _{\lambda }+\left( \lambda \left( T-t+a\right) ^{\lambda
-1}u^{2}\varphi _{\lambda }\right) _{t}+  \label{3.19}
\end{equation}%
\begin{equation*}
+\dsum\limits_{i=1}^{n}\left( -2\beta \lambda \left( T-t+a\right) ^{\lambda
-1}vv_{x_{i}}+2\lambda \left( T-t+a\right) ^{\lambda -1}vfw_{x_{i}}\right)
_{x_{i}}.
\end{equation*}%
Integrate inequality (\ref{3.19}) over $Q_{T}.$ By Gauss formula and (\ref%
{3.14}) the integral over $Q_{T}$ of the last line of (\ref{3.19}) equals
zero. Hence, we obtain the target estimate (\ref{3.121}). $\square $

\section{Lipschitz Stability and Uniqueness}

\label{sec:4}

To prove Theorem 4.1, we need in Theorems 3.1 and 3.2 the parameter $a>2$
such that $\left( T+a\right) /a^{2}<1.$ Hence, we should have 
\begin{equation*}
a>\max \left( 2,\frac{1}{2}+\sqrt{\frac{1}{4}+T}\right) .
\end{equation*}
For the sake of definiteness, we set below%
\begin{equation}
a=2+\sqrt{\frac{1}{4}+T}  \label{4.1}
\end{equation}%
and denote 
\begin{equation}
\rho =\frac{T+a}{a^{2}}=\frac{T+2+\sqrt{T+1/4}}{\left( 2+\sqrt{T+1/4}\right)
^{2}}<1.  \label{4.2}
\end{equation}

\textbf{Theorem 4.1.} \emph{Let }$N_{1},N_{2},N_{3},N_{4}>0$\emph{\ be some
numbers. Assume that in (\ref{2.1}) the function }$G=G\left( x,t,y,z\right) :%
\overline{Q}_{T}\times \mathbb{R}^{2}\rightarrow \mathbb{R}$\emph{\ is
bounded in any bounded subset of the set }$\overline{Q}_{T}\times \mathbb{R}%
^{2}$\emph{\ and such that there exist derivatives }$G_{y},G_{z}\in C\left( 
\overline{Q}_{T}\times \mathbb{R}^{2}\right) $ \emph{satisfying} 
\begin{equation}
\max \left( \sup_{Q_{T}\times \mathbb{R}^{2}}\left\vert G_{y}\left(
x,t,y,z\right) \right\vert ,\sup_{Q_{T}\times \mathbb{R}^{2}}\left\vert
G_{z}\left( x,t,y,z\right) \right\vert \right) \leq N_{1}.  \label{4.3}
\end{equation}%
\emph{Let in (\ref{2.1}) and (\ref{2.2}) the function }$K\left( x,y\right) $%
\emph{\ be bounded in }$\overline{\Omega }\times \overline{\Omega }$, \emph{%
the function }$\varkappa \in C^{1}\left( \overline{\Omega }\right) $ \emph{%
and}%
\begin{equation}
\sup_{\Omega \times \Omega }\left\vert K\left( x,y\right) \right\vert
,\left\Vert \varkappa \right\Vert _{C^{1}\left( \overline{\Omega }\right)
}\leq N_{2}.  \label{4.4}
\end{equation}%
\emph{Define sets of functions }$D_{3}\left( N_{3}\right) ,D_{4}\left(
N_{4}\right) $\emph{\ as:}%
\begin{equation}
D_{3}\left( N_{3}\right) =\left\{ u\in H_{0}^{2}\left( Q_{T}\right)
:\sup_{Q_{T}}\left\vert u\right\vert ,\sup_{Q_{T}}\left\vert \nabla
u\right\vert ,\sup_{Q_{T}}\left\vert \Delta u\right\vert \leq N_{3}\right\} ,
\label{4.5}
\end{equation}%
\begin{equation}
D_{4}\left( N_{4}\right) =\left\{ p\in H_{0}^{2}\left( Q_{T}\right)
:\sup_{Q_{T}}\left\vert p\right\vert ,\sup_{Q_{T}}\left\vert \nabla
p\right\vert \leq N_{4}\right\} ,  \label{4.6}
\end{equation}%
\emph{Let}%
\begin{equation}
N=\max \left( N_{1},N_{2},N_{3},N_{4}\right) .  \label{4.60}
\end{equation}%
\emph{Let two pairs of functions }%
\begin{equation}
\left( u_{1},p_{1}\right) ,\left( u_{2},p_{2}\right) \in D_{3}\left(
N_{3}\right) \times D_{4}\left( N_{4}\right)   \label{4.7}
\end{equation}%
\emph{\ satisfy equations (\ref{2.1}), (\ref{2.2}) as well as the following
initial and terminal conditions (see (\ref{2.4}), (\ref{2.5})):}%
\begin{equation}
u_{1}\left( x,T\right) =u_{T}^{\left( 1\right) }\left( x\right) ,\text{ }%
u_{2}\left( x,T\right) =u_{T}^{\left( 2\right) }\left( x\right) ,\text{ }%
x\in \Omega ,  \label{4.70}
\end{equation}%
\begin{equation}
u_{1}\left( x,0\right) =u_{0}^{\left( 1\right) }\left( x\right) ,\text{ }%
u_{2}\left( x,0\right) =u_{0}^{\left( 2\right) }\left( x\right) ,\text{ }%
x\in \Omega ,  \label{4.71}
\end{equation}%
\begin{equation}
p_{1}\left( x,0\right) =p_{0}^{\left( 1\right) }\left( x\right) ,\text{ }%
p_{2}\left( x,0\right) =p_{0}^{\left( 2\right) }\left( x\right) ,\text{ }%
x\in \Omega .  \label{4.72}
\end{equation}%
\emph{\ Then there exists a number }$C_{2}=C_{2}\left( \beta ,N,T\right) >0$ 
\emph{depending only on listed parameters such that the following Lipschitz
stability estimate is valid:}%
\begin{equation}
\left. 
\begin{array}{c}
\left\Vert u_{1}-u_{2}\right\Vert _{H^{1,0}\left( Q_{T}\right) }+\left\Vert
p_{1}-p_{2}\right\Vert _{H^{1,0}\left( Q_{T}\right) }\leq  \\ 
\leq C_{2}\left( \left\Vert u_{T}^{\left( 1\right) }-u_{T}^{\left( 2\right)
}\right\Vert _{H^{1}\left( \Omega \right) }+\left\Vert u_{0}^{\left(
1\right) }-u_{0}^{\left( 2\right) }\right\Vert _{L_{2}\left( \Omega \right)
}+\left\Vert p_{0}^{\left( 1\right) }-p_{0}^{\left( 2\right) }\right\Vert
_{L_{2}\left( \Omega \right) }\right) .%
\end{array}%
\right.   \label{4.11}
\end{equation}%
\emph{Next, if in (\ref{4.70})-(\ref{4.72})} 
\begin{equation*}
u_{T}^{\left( 1\right) }\left( x\right) =u_{T}^{\left( 2\right) }\left(
x\right) ,\text{ }u_{0}^{\left( 1\right) }\left( x\right) =u_{0}^{\left(
2\right) }\left( x\right) ,\text{ }p_{0}^{\left( 1\right) }\left( x\right)
=p_{0}^{\left( 2\right) }\left( x\right) ,\text{ }x\in \Omega ,
\end{equation*}%
\emph{then }$u_{1}\left( x,t\right) \equiv u_{2}\left( x,t\right) $\emph{\
and }$p_{1}\left( x,t\right) \equiv p_{2}\left( x,t\right) $ \emph{in }$%
Q_{T},$\emph{\ which means that problem (\ref{2.1})-(\ref{2.5}) has at most
one solution satisfying conditions of this theorem.}

\textbf{Remark 4.1.}\emph{\ Imposing condition (\ref{4.7}) that functions }$%
u_{i},p_{i},i=1,2$\emph{\ belong to  a priori prescribed bounded sets is
typical in the field of Ill-Posed and Inverse Problems. This is due to the
regularization theory, see, e.g. \cite{T}. }

\textbf{Proof}. In this proof the parameter $a$ in the function $\varphi
_{\lambda }\left( t\right) =e^{2\left( T-t+a\right) ^{\lambda }}$ is the
same as in (\ref{4.1}). Below $C_{2}=C_{2}\left( \beta ,N,T\right) >0$
denotes different numbers depending only on listed parameters. For any four
numbers $a_{1},b_{1},a_{2},b_{2}\in \mathbb{R},$ let $\widetilde{a}%
=a_{1}-a_{2}$ and $\widetilde{b}=b_{1}-b_{2}.$ Then 
\begin{equation}
a_{1}b_{1}-a_{2}b_{2}=\widetilde{a}b_{1}+\widetilde{b}a_{2}.  \label{4.12}
\end{equation}

Denote 
\begin{equation}
\widetilde{u}\left( x,t\right) =u_{1}\left( x,t\right) -u_{2}\left(
x,t\right) ,\text{ }\widetilde{p}\left( x,t\right) =p_{1}\left( x,t\right)
-p_{2}\left( x,t\right) ,\text{ }\left( x,t\right) \in Q_{T},  \label{4.13}
\end{equation}%
\begin{equation}
\widetilde{u}_{T}\left( x\right) =u_{T}^{\left( 1\right) }\left( x\right)
-u_{T}^{\left( 2\right) }\left( x\right) ,\text{ }\widetilde{u}_{0}\left(
x\right) =u_{0}^{\left( 1\right) }\left( x\right) -u_{0}^{\left( 2\right)
}\left( x\right) ,\text{ }x\in \Omega ,  \label{4.14}
\end{equation}%
\begin{equation}
\widetilde{p}_{0}\left( x\right) =p_{0}^{\left( 1\right) }\left( x\right)
-p_{0}^{\left( 2\right) }\left( x\right) ,\text{ }x\in \Omega .  \label{4.15}
\end{equation}

Using the multidimensional analog of Taylor formula \cite{V}, (\ref{4.3})
and (\ref{4.13}), we obtain%
\begin{equation}
\left. 
\begin{array}{c}
G\left( x,t,\dint\limits_{\Omega }K\left( x,y\right) p_{1}\left( y,t\right)
dy,p_{1}\left( x,t\right) \right) - \\ 
-G\left( x,t,\dint\limits_{\Omega }K\left( x,y\right) p_{2}\left( y,t\right)
dy,p_{2}\left( x,t\right) \right) = \\ 
=F_{1}\left( x,t\right) \dint\limits_{\Omega }K\left( x,y\right) \widetilde{p%
}\left( y,t\right) dy+F_{2}\left( x,t\right) \widetilde{p}\left( x,t\right) ,%
\end{array}%
\right.  \label{4.16}
\end{equation}%
where functions $F_{1}$,$F_{2}\in C\left( \overline{Q}_{T}\right) $ and by (%
\ref{4.60}) 
\begin{equation}
\left\vert F_{1}\left( x,t\right) \right\vert ,\left\vert F_{2}\left(
x,t\right) \right\vert \leq N,\text{ }\left( x,t\right) \in \overline{Q}_{T}.
\label{4.17}
\end{equation}

Subtract equations (\ref{2.1}), (\ref{2.2}) for the pair $\left(
u_{2},p_{2}\right) $ from the same equations for the pair $\left(
u_{1},p_{1}\right) .$ Using (\ref{4.4})-(\ref{4.7}), (\ref{4.12})-(\ref{4.17}%
) and turning resulting two equations in two inequalities, we obtain for $%
\left( x,t\right) \in Q_{T}$ 
\begin{equation}
\left\vert \widetilde{u}_{t}+\beta \Delta \widetilde{u}\right\vert \left(
x,t\right) \leq C_{2}\left( \left\vert \nabla \widetilde{u}\left( x,t\right)
\right\vert +\left\vert \widetilde{p}\left( x,t\right) \right\vert
+\dint\limits_{\Omega }\left\vert \widetilde{p}\left( y,t\right) \right\vert
dy\right) ,  \label{4.18}
\end{equation}%
\begin{equation}
\left\vert \widetilde{p}_{t}-\beta \Delta \widetilde{p}+\varkappa
^{2}p_{1}\Delta \widetilde{u}\right\vert \left( x,t\right) \leq C_{2}\left(
\left\vert \nabla \widetilde{p}\left( x,t\right) \right\vert +\left\vert 
\widetilde{p}\left( x,t\right) \right\vert +\left\vert \nabla \widetilde{u}%
\left( x,t\right) \right\vert \right) .  \label{4.19}
\end{equation}%
Initial and terminal conditions for functions $\widetilde{u}$ and $%
\widetilde{p}$ are:%
\begin{equation}
\widetilde{u}\left( x,T\right) =\widetilde{u}_{T}\left( x\right) ,\widetilde{%
u}\left( x,0\right) =\widetilde{u}_{0}\left( x\right) ,\text{ }\widetilde{p}%
\left( x,0\right) =\widetilde{p}_{0}\left( x\right) .  \label{4.20}
\end{equation}%
In addition, functions 
\begin{equation}
\widetilde{u},\widetilde{p}\in H_{0}^{2}\left( Q_{T}\right) .  \label{4.200}
\end{equation}

Squaring both sides of each of inequalities (\ref{4.18}) and (\ref{4.19}),
multiplying them by the unction $\varphi _{\lambda }\left( t\right) ,$ using
Cauchy-Schwarz inequality and integrating over $Q_{T},$ we obtain:%
\begin{equation}
\dint\limits_{Q_{T}}\left( \widetilde{u}_{t}+\beta \Delta \widetilde{u}%
\right) ^{2}\varphi _{\lambda }dxdt\leq C_{2}\dint\limits_{Q_{T}}\left(
\left( \nabla \widetilde{u}\right) ^{2}+\widetilde{p}^{2}+\dint\limits_{%
\Omega }\widetilde{p}^{2}\left( y,t\right) dy\right) \varphi _{\lambda }dxdt,
\label{4.21}
\end{equation}%
and also 
\begin{equation*}
\dint\limits_{Q_{T}}\left( \widetilde{p}_{t}-\beta \Delta \widetilde{p}%
+\varkappa ^{2}p_{1}\Delta \widetilde{u}\right) ^{2}\varphi _{\lambda
}dxdt\leq 
\end{equation*}%
\begin{equation}
\leq C_{2}\dint\limits_{Q_{T}}\left( \left( \nabla \widetilde{p}\right) ^{2}+%
\widetilde{p}^{2}+\left( \nabla \widetilde{u}\right) ^{2}\right) \varphi
_{\lambda }dxdt.  \label{4.22}
\end{equation}

It follows from (\ref{4.200}) that we can apply the Carleman estimate of
Theorem 3.1 to the left hand side of (\ref{4.21}) and the quasi Carleman
estimate of Theorem 3.2 to the left hand side of (\ref{4.22}). First, we
work with (\ref{4.21}). Using  (\ref{3.300}) and (\ref{4.21}), we obtain%
\begin{equation*}
\frac{2}{3}\sqrt{\lambda }\beta \dint\limits_{Q_{T}}\left( \nabla \widetilde{%
u}\right) ^{2}\varphi _{\lambda }dxdt+\frac{\lambda ^{2}}{12}a^{\lambda
-2}\dint\limits_{Q_{T}}\widetilde{u}^{2}\varphi _{\lambda }dxdt-
\end{equation*}%
\begin{equation}
-\frac{2}{3}e^{2a^{\lambda }}\dint\limits_{\Omega }\left( \beta \left(
\nabla \widetilde{u}\right) ^{2}+\frac{1}{2}\widetilde{u}^{2}\right) \left(
x,T\right) -\frac{2}{3}\lambda \left( T+a\right) ^{\lambda
-1}\dint\limits_{\Omega }\widetilde{u}^{2}\left( x,0\right) dx\leq 
\label{4.23}
\end{equation}%
\begin{equation*}
\leq C_{2}\dint\limits_{Q_{T}}\left( \left( \nabla \widetilde{u}\right) ^{2}+%
\widetilde{p}^{2}+\dint\limits_{\Omega }\widetilde{p}^{2}\left( y,t\right)
dy\right) \varphi _{\lambda }dxdt,\text{ }\forall \lambda \geq \lambda _{0}.
\end{equation*}%
Let $\lambda _{0}>64$\ be the number defined in (\ref{3.9}). Choose a
sufficiently large number $\lambda _{1}=\lambda _{1}\left( \beta ,N,T\right)
\geq \lambda _{0}$ such that 
\begin{equation*}
\frac{1}{3}\sqrt{\lambda _{1}}\beta \geq C_{2}.
\end{equation*}%
Then (\ref{4.23}) leads to:%
\begin{equation*}
\frac{1}{3}\sqrt{\lambda }\beta \dint\limits_{Q_{T}}\left( \nabla \widetilde{%
u}\right) ^{2}\varphi _{\lambda }dxdt+\frac{\lambda ^{2}}{12}a^{\lambda
-2}\dint\limits_{Q_{T}}\widetilde{u}^{2}\varphi _{\lambda }dxdt\leq 
\end{equation*}%
\begin{equation}
\leq C_{2}\dint\limits_{Q_{T}}\widetilde{p}^{2}\varphi _{\lambda }dxdt+
\label{4.24}
\end{equation}%
\begin{equation*}
+C_{2}e^{2a^{\lambda }}\dint\limits_{\Omega }\left( \left( \nabla \widetilde{%
u}\right) ^{2}+\widetilde{u}^{2}\right) \left( x,T\right) +C_{2}\lambda
\left( T+a\right) ^{\lambda -1}\dint\limits_{\Omega }\widetilde{u}^{2}\left(
x,0\right) dx,\text{ }\forall \lambda \geq \lambda _{1}.
\end{equation*}

To handle terms with the function $\widetilde{p}$ in the second line of (\ref%
{4.24}), we apply now quasi-Carleman estimate (\ref{3.121}) to the first
line of (\ref{4.22}). We set in (\ref{3.121}):%
\begin{equation}
u=\widetilde{p},q=\widetilde{u},f=\varkappa ^{2}p_{1}.  \label{4.240}
\end{equation}%
Keeping in mind (\ref{4.4})-(\ref{4.7}), we obtain from (\ref{3.121}), (\ref%
{4.22}) and (\ref{4.240}):   
\begin{equation*}
\frac{\lambda ^{2}}{4}a^{2\lambda -2}\dint\limits_{Q_{T}}\widetilde{p}%
^{2}\varphi _{\lambda }dxdt+\beta \lambda a^{\lambda
-1}\dint\limits_{Q_{T}}\left( \nabla \widetilde{p}\right) ^{2}\varphi
_{\lambda }dxdt\leq 
\end{equation*}%
\begin{equation}
\leq C_{2}\lambda \left( T+a\right) ^{\lambda -1}\dint\limits_{Q_{T}}\left(
\nabla \widetilde{u}\right) ^{2}\varphi _{\lambda }dxdt+  \label{4.25}
\end{equation}%
\begin{equation*}
+C_{2}\dint\limits_{Q_{T}}\left( \left( \nabla \widetilde{p}\right) ^{2}+%
\widetilde{p}^{2}+\left( \nabla \widetilde{u}\right) ^{2}\right) \varphi
_{\lambda }dxdt-
\end{equation*}%
\begin{equation*}
+\lambda \left( T+a\right) ^{\lambda -1}e^{2\left( T+a\right) ^{\lambda
}}\dint\limits_{\Omega }\widetilde{p}^{2}\left( x,0\right) dx,\text{ }%
\forall \lambda \geq \lambda _{0}.
\end{equation*}%
By (\ref{4.1}) we can choose $\lambda _{2}=\lambda _{2}\left( \beta
,N,T\right) \geq \lambda _{1}$ such that%
\begin{equation*}
\min \left( \frac{\lambda ^{2}}{8}a^{2\lambda -2},\frac{\beta \lambda
a^{\lambda -1}}{2}\right) \geq C_{2},\text{ }\forall \lambda \geq \lambda
_{2}.
\end{equation*}%
Hence, (\ref{4.25}) implies%
\begin{equation*}
\lambda ^{2}a^{2\lambda -2}\dint\limits_{Q_{T}}\widetilde{p}^{2}\varphi
_{\lambda }dxdt+\lambda a^{\lambda -1}\dint\limits_{Q_{T}}\left( \nabla 
\widetilde{p}\right) ^{2}\varphi _{\lambda }dxdt\leq 
\end{equation*}%
\begin{equation}
\leq C_{2}\lambda \left( T+a\right) ^{\lambda -1}\dint\limits_{Q_{T}}\left(
\nabla \widetilde{u}\right) ^{2}\varphi _{\lambda }dxdt+  \label{4.26}
\end{equation}%
\begin{equation*}
+C_{2}\lambda \left( T+a\right) ^{\lambda -1}e^{2\left( T+a\right) ^{\lambda
}}\dint\limits_{\Omega }\widetilde{p}^{2}\left( x,0\right) dx,\text{ }%
\forall \lambda \geq \lambda _{2}.
\end{equation*}%
In particular, (\ref{4.1}), (\ref{4.2}) and (\ref{4.26}) imply%
\begin{equation}
\dint\limits_{Q_{T}}\widetilde{p}^{2}\varphi _{\lambda }dxdt\leq C_{2}\frac{%
\rho ^{\lambda -1}}{\lambda }\dint\limits_{Q_{T}}\left( \nabla \widetilde{u}%
\right) ^{2}\varphi _{\lambda }dxdt+  \label{4.27}
\end{equation}%
\begin{equation*}
+C_{2}\frac{\rho ^{\lambda -1}}{\lambda }e^{2\left( T+a\right) ^{\lambda
}}\dint\limits_{\Omega }\widetilde{p}^{2}\left( x,0\right) dx,\text{ }%
\forall \lambda \geq \lambda _{2}.
\end{equation*}%
Comparing (\ref{4.27}) with the second line of (\ref{4.24}), we obtain%
\begin{equation*}
\frac{1}{3}\sqrt{\lambda }\beta \dint\limits_{Q_{T}}\left( \nabla \widetilde{%
u}\right) ^{2}\varphi _{\lambda }dxdt+\frac{\lambda ^{2}}{12}a^{\lambda
-2}\dint\limits_{Q_{T}}\widetilde{u}^{2}\varphi _{\lambda }dxdt\leq 
\end{equation*}%
\begin{equation}
\leq C_{2}\frac{\rho ^{\lambda -1}}{\lambda }\dint\limits_{Q_{T}}\left(
\nabla \widetilde{u}\right) ^{2}\varphi _{\lambda }dxdt+  \label{4.28}
\end{equation}%
\begin{equation*}
+C_{2}\frac{\rho ^{\lambda -1}}{\lambda }e^{2\left( T+a\right) ^{\lambda
}}\dint\limits_{\Omega }\widetilde{p}^{2}\left( x,0\right) dx+
\end{equation*}%
\begin{equation*}
+C_{2}e^{2a^{\lambda }}\dint\limits_{\Omega }\left( \left( \nabla \widetilde{%
u}\right) ^{2}+\widetilde{u}^{2}\right) \left( x,T\right) +C_{2}\lambda
\left( T+a\right) ^{\lambda -1}\dint\limits_{\Omega }\widetilde{u}^{2}\left(
x,0\right) dx,\text{ }\forall \lambda \geq \lambda _{2}.
\end{equation*}%
Since $\rho <1$ by (\ref{4.2}), then we can choose $\lambda _{3}=\lambda
_{3}\left( \beta ,N,T\right) \geq \lambda _{2}$ such that 
\begin{equation*}
\frac{1}{6}\sqrt{\lambda }\beta \geq C_{2}\frac{\rho ^{\lambda -1}}{\lambda }%
,\text{ }\forall \lambda \geq \lambda _{3}.
\end{equation*}%
Hence, (\ref{4.28}) implies%
\begin{equation}
\dint\limits_{Q_{T}}\left( \left( \nabla \widetilde{u}\right) ^{2}+%
\widetilde{u}^{2}\right) \varphi _{\lambda }dxdt\leq C_{2}e^{2\left(
T+a\right) ^{\lambda }}\dint\limits_{\Omega }\widetilde{p}^{2}\left(
x,0\right) dx+  \label{4.29}
\end{equation}%
\begin{equation*}
+C_{2}e^{2a^{\lambda }}\dint\limits_{\Omega }\left( \left( \nabla \widetilde{%
u}\right) ^{2}+\frac{1}{2}\widetilde{u}^{2}\right) \left( x,T\right)
+C_{2}\lambda \left( T+a\right) ^{\lambda -1}\dint\limits_{\Omega }%
\widetilde{u}^{2}\left( x,0\right) dx,\text{ }\forall \lambda \geq \lambda
_{3}.
\end{equation*}%
Comparing (\ref{4.29}) with (\ref{4.26}), we incorporate in estimate (\ref%
{4.29}) the term with $\left( \nabla \widetilde{p}\right) ^{2}+\widetilde{p}%
^{2},$ 
\begin{equation*}
\dint\limits_{Q_{T}}\left( \left( \nabla \widetilde{u}\right) ^{2}+%
\widetilde{u}^{2}+\left( \nabla \widetilde{p}\right) ^{2}+\widetilde{p}%
^{2}\right) \varphi _{\lambda }dxdt\leq 
\end{equation*}%
\begin{equation}
\leq C_{2}\left( 1+\frac{T}{a}\right) ^{\lambda -1}e^{2\left( T+a\right)
^{\lambda }}\dint\limits_{\Omega }\widetilde{p}^{2}\left( x,0\right) dx+
\label{4.30}
\end{equation}%
\begin{equation*}
+C_{2}\left( 1+\frac{T}{a}\right) ^{\lambda -1}e^{2a^{\lambda
}}\dint\limits_{\Omega }\left( \left( \nabla \widetilde{u}\right) ^{2}+%
\widetilde{u}^{2}\right) \left( x,T\right) +
\end{equation*}%
\begin{equation*}
+C_{2}\lambda \left( T+a\right) ^{2\lambda -2}\dint\limits_{\Omega }%
\widetilde{u}^{2}\left( x,0\right) dx,\text{ }\forall \lambda \geq \lambda
_{3}.
\end{equation*}%
By (\ref{3.2})%
\begin{equation*}
\dint\limits_{Q_{T}}\left( \left( \nabla \widetilde{u}\right) ^{2}+%
\widetilde{u}^{2}+\left( \nabla \widetilde{p}\right) ^{2}+\widetilde{p}%
^{2}\right) \varphi _{\lambda }dxdt\geq 
\end{equation*}%
\begin{equation*}
\geq e^{2a^{\lambda }}\dint\limits_{Q_{T}}\left( \left( \nabla \widetilde{u}%
\right) ^{2}+\widetilde{u}^{2}+\left( \nabla \widetilde{p}\right) ^{2}+%
\widetilde{p}^{2}\right) dxdt.
\end{equation*}%
Hence, dividing (\ref{4.30}) by $e^{2a^{\lambda }},$ setting in the
resulting inequality $\lambda =\lambda _{3}=\lambda _{3}\left( \beta
,N,T\right) $ and using (\ref{4.20}), we obtain%
\begin{equation}
\left\Vert \widetilde{u}\right\Vert _{H^{1,0}\left( Q_{T}\right)
}^{2}+\left\Vert \widetilde{p}\right\Vert _{H^{1,0}\left( Q_{T}\right)
}^{2}\leq C_{2}\left( \left\Vert \widetilde{u}_{T}\right\Vert _{H^{1}\left(
\Omega \right) }^{2}+\left\Vert \widetilde{u}_{0}\right\Vert _{L_{2}\left(
\Omega \right) }^{2}+\left\Vert \widetilde{p}_{0}\right\Vert _{L_{2}\left(
\Omega \right) }^{2}\right) .  \label{4.31}
\end{equation}%
The target estimate (\ref{4.11}) of this theorem follows immediately from a
combination of (\ref{4.31}) with (\ref{4.13})-(\ref{4.15}). $\square $

\end{document}